\documentclass[11pt]{article}
\usepackage{amssymb,latexsym}
\newtheorem{theo}{Theorem}

\newtheorem{lemm}[theo]{Lemma}
\newtheorem{coro}[theo]{Corollary}

\def\nn{\nonumber}
\def\ds{\displaystyle}
\def\C{\mathbb{C}} 
\def\Z{\mathbb{Z}} 
\def\ra{\rangle}
\def\la{\langle}
\def\lb{[\![}
\def\rb{]\!]}
\def\l{\ldots}
\def\t{\theta}
\def\q{{\bar q}}
\def\L{{\bar L}}

\def\mybox{\hfill$\Box$}
\begin{document}
\addtolength{\baselineskip}{2mm}
\addtolength{\abovedisplayskip}{1mm}
\addtolength{\belowdisplayskip}{1mm}
\addtolength{\parskip}{2mm}
\begin{center} 
{\Large \bf 
Jacobson generators of the quantum superalgebra\\[1mm]
$U_q[sl(n+1|m)]$ and Fock representations}\\[5mm] 
{\bf T.D.\ Palev}\footnote{E-mail~: tpalev@inrne.bas.bg.}\\[1mm]
Institute for Nuclear Research and Nuclear Energy,\\ 
Boul.\ Tsarigradsko Chaussee 72, 1784 Sofia, Bulgaria;\\[2mm] 
{\bf N.I.\ Stoilova}\footnote{E-mail~: ptns@pt.tu-clausthal.de.
Permanent address~: Institute for Nuclear Research and Nuclear Energy, 
Boul.\ Tsarigradsko Chaussee 72, 1784 Sofia, Bulgaria.}\\[1mm]
Mathematical Physics Group, Department of Physics, 
Technical University of Clausthal,\\
Leibnizstrasse~10, D-38678 Clausthal-Zellerfeld, Germany;\\[2mm] 
{\bf and J.\ Van der Jeugt}\footnote{E-mail~:
Joris.VanderJeugt@rug.ac.be.}\\[1mm]
Department of Applied Mathematics and Computer Science, 
University of Ghent,\\
Krijgslaan 281-S9, B-9000 Gent, Belgium. 
\end{center}

\vskip 10mm
\noindent Corresponding author~: J.\ Van der Jeugt. 
Tel~: ++32 9 2644812. Fax~: ++32 9 2644995. 
E-mail~: Joris.VanderJeugt@rug.ac.be
\vskip 3mm
\noindent Running title~: Jacobson generators of $U_q[sl(n+1|m)]$

\vskip 10mm
\begin{abstract}
As an alternative to Chevalley generators,
we introduce Jacobson generators for the quantum superalgebra
$U_q[sl(n+1|m)]$. The expressions of all Cartan-Weyl elements
of $U_q[sl(n+1|m)]$ in terms of these Jacobson generators 
become very simple. 
We determine and prove certain triple relations between
the Jacobson generators, necessary for a complete set of
supercommutation relations between the Cartan-Weyl elements.
Fock representations are defined, and a substantial part of this
paper is devoted to the computation of the action of Jacobson generators
on basis vectors of these Fock spaces. 
It is also determined when these Fock representations are unitary.
Finally, Dyson and Holstein-Primakoff realizations are given,
not only for the Jacobson generators, but for all Cartan-Weyl
elements of $U_q[sl(n+1|m)]$.
\end{abstract} 
\vskip 2mm

\renewcommand{\thesection}{\Roman{section}}

\section{Introduction}
\setcounter{equation}{0}

The quantization of simple Lie algebras~\cite{Drinfeld,Jimbo} or Lie
superalgebras~\cite{Bracken,Chaichian,Floreanini1,Khoroshkin,Scheunert}
as quasitriangular Hopf (super)algebras has been carried out
more than a decade ago.
Since then, these structures have received much attention both
in the mathematical and physical literature.
In a physical context, one is mostly dealing with 
representations or realizations of these quantized algebras.
This is in fact the main topic of the present paper~: 
certain special representations (Fock representations)
and related realizations (Dyson and Holstein-Primakoff)
of the quantum superalgebra $U_q[sl(n+1|m)]$ are presented.

The Lie superalgebra $sl(n+1|m)$ is one of the basic classical
simple Lie superalgebras in Kac's classification~\cite{Kac1}.
It can be considered as the superanalogue of the special linear
Lie algebra $sl(n+1)$. 
The quantum superalgebra $U_q[sl(n+1|m)]$ is a Hopf superalgebra
deformation of the associative superalgebra $U[sl(n+1|m)]$,
the universal enveloping superalgebra of $sl(n+1|m)$.
At this point, it is already worth observing that the more familiar
case of $sl(n+1)$ and $U_q[sl(n+1)]$ just follows by putting
$m=0$. The readers who are interested in this case only can
still use all formulas presented in this paper, simply taking $m$
equal to 0.

For a definition of the quantum superalgebra $U_q[sl(n+1|m)]$,
we refer to~\cite{Bracken,Chaichian,Floreanini1,Khoroshkin,Scheunert}.
Usually, $U_q[sl(n+1|m)]$ is defined by its Chevalley generators
(often denoted by $e_i$, $f_i$ and $h_i$, with $i=1,\ldots,n+m$),
subject to the Cartan-Kac relations and the Serre 
relations~\cite{Floreanini1,Khoroshkin,Scheunert}.
Besides these defining relations, also the other Hopf
superalgebra maps (comultiplication, co-unit and antipode) are
part of the definition. In this paper, however, we do not use
these other Hopf superalgebra maps; 
so we shall concentrate on $U_q[sl(n+1|m)]$ as an
associative superalgebra.

The definition in terms of Chevalley generators has the advantage
that the comultiplication, co-unit and antipode are easy to give.
Furthermore, certain representations can be constructed
explicitly (e.g.\ for the essentially typical representations 
a Gelfand-Zetlin basis exist for which the action of the Chevalley
generators is known~\cite{PSV94}).
Having certain physical applications in mind, however, it is
sometimes more useful to work with a different set of generators
for $U_q[sl(n+1|m)]$.

The different set of generators for $U_q[sl(n+1|m)]$ given here
are the Jacobson generators (denoted by $a_i^+$, $a_i^-$ and $H_i$,
with $i=1,\ldots,n+m$). For the case of $sl(n+1)$, such generators
were originally introduced by Jacobson~\cite{Jacobson,PV00}.
The use of Jacobson generators has a number of advantages.

First of all, in certain applications it is necessary to have a
complete basis of $U_q[sl(n+1|m)]$ (following from the 
Poincar\'e-Birkhoff-Witt theorem). Such a basis is given in
terms of the Cartan-Weyl elements. Although it is possible to
express all Cartan-Weyl elements in terms of the Chevalley
generators, such expressions soon become rather unmanageable.
In terms of the Jacobson generators, the description of
all Cartan-Weyl elements is very easy.

Secondly, the Jacobson generators allow for the construction of a class
of irreducible $U_q[sl(n+1|m)]$ modules $W_p$, $p \in \C$,
called Fock representations. 
The Fock representations corresponding to different
$p$ are inequivalent. For $p$ a positive integer
they provide an explicit construction (basis and transformation 
of the basis under the action of the generators) 
of (deformations of) atypical
representations of $U_q[sl(n+1|m)]$. This is an interesting mathematical
result, since even in the nondeformed case all atypical representations of
$sl(n+1|m)$ were not explicitly
constructed so far (e.g.\ even a dimension formula is unknown).

A disadvantage of the Jacobson generators compared to the
Chevalley generators is that the explicit expressions for 
the other Hopf (super)algebra maps
(comultiplication, co-unit and antipode) become very
lengthy and complicated.

The results of the present paper provide a mathematical background for
further studies of noncanonical quantum statistics initiated 
in~\cite{Palev-thesis}  (see also~\cite{PV00} and~\cite{JPV} 
for further references). The approach is
based on the concept of creation and annihilation operators
(CAO's) of a simple Lie (super)algebra ${\cal A}$ and its
Fock representations~\cite{Palev80}. The CAO's of ${\cal A}$
provide a description of ${\cal A}$ in terms of generators and relations,
which are different from the Chevalley generators. 
In this terminology any $n$ pairs of para-Fermi operators~\cite{Green}  
are CAO's of $so(2n+1)$~\cite{KR}  
and any $n$ pairs of para-Bose operators~\cite{Green} 
are CAO's of the orthosymplectic Lie superalgebra  
$osp(1|2n)$~\cite{Ganchev}. The CAO's of $sl(n+1)$~\cite{Palev-thesis} 
and of $sl(1|n)$~\cite{Palev80} lead to new quantum statistics.
Generalizing the results of Jacobson on Lie triple systems~\cite{Jacobson},
Okubo has reformulated all above statistics in terms of
Lie supertriple systems~\cite{Okubo}. 
In this setting the CAO's of the Lie (super)algebras mentioned above are
generators of the related (super)triple systems. 
This is another reason to call them  Jacobson generators (JG's). 
The link between the JG's and the simple Lie superalgebras
provides a natural background for their $q$-deformations
(we refer to~\cite{Palev-ospnm} for more discussion in this respect).

The representations of (quantum) superalgebras have certainly wider
applications. These algebras
(and in particular $U_q[gl(n+1|m)]$~\cite{Mark}) 
play a role for finding new solutions of the quantum Yang-Baxter equations
and for the construction of solvable models. As examples we mention the
supersymmetric solvable $t-J$ models of correlated electrons~\cite{Sarker} 
and their quantum analogue~\cite{Fei}. Some other
potential physical applications are mentioned in the last section.

In section~II we define the Jacobson generators of $U_q[sl(n+1|m)]$,
as a special subset of the Cartan-Weyl elements.
The description of all Cartan-Weyl elements in terms of the
Jacobson generators becomes very simple (Theorem~1). 
However, in order to apply these results (e.g.\ in representations)
one must have a list of all (super)commutation relations
between these Cartan-Weyl elements; in terms of Jacobson
generators, this means one has to determine certain triple
relations. These are also given in Theorem~1, together with
their proof. 

In section~III we define Fock representations for $U_q[sl(n+1|m)]$,
related to the earlier defined Jacobson generators.
The main part of this section is devoted to the proof of Theorem~2,
describing the action of the Jacobson generators on a basis
of the Fock representation. This proof is rather technical
and lengthy, and has been divided in a number of lemmas. 
The essential result is that these Fock representations
are labelled by a number $p$; when $p$ is a nonnegative
integer, the Fock representation is finite-dimensional.
%For an application of these Fock representations, see~\cite{PSV00}~:
%herein, the Jacobson generators are interpreted as operators
%creating or annihilating a ``particle'', and the underlying
%quantum statistics is discussed.

The Fock representations determined in section~III are further 
analysed in section~IV. In particular, following conditions
required in a physical context, it is determined when these
Fock representations are unitary (or unitarizable, or Hermitian),
see Theorem~6. In that case, an orthonormal basis of the Fock
space is given, together with the action of the Jacobson 
generators on these basis elements.

Inspired by the Fock representations, we can give new 
expressions for the Dyson
and Holstein-Primakoff realizations of $U_q[sl(n+1|m)]$ (section~V).
In~\cite{Palev99}, the Dyson and Holstein-Primakoff realizations
for the Chevalley generators of $U_q[gl(n|m)]$ was already
given. Here, we give Dyson and Holstein-Primakoff realizations
for the Jacobson generators of $U_q[sl(n+1|m)]$ (Theorems~7 and~8);
from these, the corresponding realization for all Cartan-Weyl
elements are deduced. 
All these realizations are in terms of $n$ pairs of Bose and $m$
pairs of Fermi creation and annihilation operators.
The Holstein-Primakoff realization becomes particularly
simple when expressed in terms of $q$-deformed Bose and
Fermi creation and annihilation operators.

Unless otherwise stated, we consider in this paper $U_q[sl(n+1|m)]$
as a module over the algebra $\C[[h]]$ (with $q=e^h$)
of formal power series over an indeterminate $h$.
It is important to note however that all considerations remain
true if one replaces $h$ by a complex number such that $q=e^h$ is not
a root of unity.
In fact, most of our results hold also for $q$ being a root of~1,
including the unitary Fock representations and the Dyson and
Holstein-Primakoff realizations.

Throughout the paper we use the notation~: 
JGs for Jacobson generators;
$\Z$ (resp.\ $\Z_+$) for the set of all integers (resp.\ of all
nonnegative integers);
$\Z_2=\{\bar{0},\bar{1}\}$ for the ring of all integers modulo 2;
$\C$ for all complex numbers. 
Furthermore~:
\begin{eqnarray}
&& [x]={q^x-q^{-x} \over q-q^{-1}}, \hbox{ when }x\in\C ,\\ 
&&[r;s]=\{r,r+1,r+2,\l,s-1,s\}, \hbox{ for } r\le s\in \Z ; \\
&&\t_i= \left\{ \begin{array}{lll}
 {\bar 0} & \hbox{if} & i\in [0;n] \\ 
 {\bar 1} & \hbox{if} & i\in[n+1;n+m]
 \end{array}\right. ; \quad
 \t_{ij}=\t_i+\t_j; \\
&&[a,b]=ab-ba,\;\; \{a,b\}=ab+ba,
\;\;\lb a,b\rb=ab-(-1)^{\deg(a)\deg(b)}ba; \\
&&[a,b]_x=ab-xba,\;\; \{a,b\}_x=ab+xba,
\;\;\lb a,b\rb_x=ab-(-1)^{\deg(a)\deg(b)}xba, 
\end{eqnarray}
where $\deg(a)\in\Z_2$ refers to the degree (or grading) of $a$ when
$a$ is a homogeneous element of a superalgebra.

\section{Jacobson generators of $U_q[sl(n+1|m)]$}
\setcounter{equation}{0}

Although the quantization ($q$-deformation) of simple Lie algebras
and basic Lie superalgebras is usually carried out
in terms of their Chevalley generators,
there exist alternative descriptions in terms
of so-called deformed creation and annihilation
operators for the $q$-deformation of $osp(1|2n)$~\cite{Palev-osp12n}, 
$so(2n+1)$~\cite{Palev-so}, $osp(2n+1|m)$~\cite{Palev-ospnm}, 
$sl(n+1)$~\cite{PP} and $sl(n+1|m)$~\cite{PS99}.
These alternative generators have the advantage that
in some natural interpretation they have a direct physical significance;
furthermore, they allow the definition and construction
of a mathematically interesting and physically important
class of irreducible representations, the Fock representations.

The Hopf superalgebra $U_q[sl(n+1|m)]$ is defined in the sense
of Drinfeld~\cite{Drinfeld}, as a topologically free $\C[[h]]$ module.
As a superalgebra, $U_q[sl(n+1|m)]$ is usually 
defined by means of its Chevalley
generators, subject to the Cartan-Kac relations and the Serre
relations~\cite{Floreanini1,Khoroshkin,Scheunert}.
Unlike the Lie algebra case, there is an ``extra Serre relation''
involving the generator associated with an odd simple 
root~\cite{Floreanini1,Khoroshkin,Scheunert,Scheunert2}.
This property was investigated further by Yamane~\cite{Yamane1,Yamane2}.
Indeed, for the basic classical Lie superalgebras there exist
many non-isomorphic simple root systems~\cite{Kac1};
one of these, having only one odd simple root, is known
as the distinguished simple root system~\cite{Kac1}.
The classical description of $U_q[sl(n+1|m)]$ is in terms of relations
and generators associated with this distinguished simple root system.
Yamane~\cite{Yamane1,Yamane2} studied Hopf superalgebras
in terms of relations and generators associated with other
simple root systems. Apparently, this gives rise to more involved
extra Serre relations. Moreover, the structure of the Hopf
superalgebra seems to depend on the choice of simple root 
system~\cite{Yamane1,Yamane2}.
In this paper, $U_q[sl(n+1|m)]$ stands for the usual Hopf superalgebra
associated with the distinguished simple root system.
But we shall be dealing with an alternative
set of generators (and relations) for $U_q[sl(n+1|m)]$.

In this section we shall recall the definition of
deformed creation and annihilation operators of $U_q[sl(n+1|m)]$,
and refer to them as Jacobson generators (JGs) since they are 
closely related to generators in the sense of a Lie
supertriple system~\cite{Okubo} (and for Lie triple systems, such
generators were originally introduced by Jacobson~\cite{Jacobson,PV00}).
The definition of JGs can be best presented in the framework
of a set of Cartan-Weyl elements of $U_q[sl(n+1|m)]$.
Furthermore, in order to construct the Fock representations
explicitly, it is necessary to have a complete list of 
so-called triple relations between the JGs.
Such relations can be deduced from the supercommutation
relations between all Cartan-Weyl elements. 
So we begin this section by recalling some properties
of Cartan-Weyl elements of $U_q[gl(n+1|m)]$, 
deduced in~\cite{PalevTolstoy}, which are then easily restricted
to the case of $U_q[sl(n+1|m)]$.

Although a set of generators, such as the Chevalley generators,
is sufficient for the definition of $U_q[gl(n+1|m]$ as an
associative algebra, it is not sufficient for describing a 
basis of $U_q[gl(n+1|m]$. For this purpose, the 
construction of a set of Cartan-Weyl elements is necessary.
For $U_q[gl(n+1|m]$, a set of Cartan-Weyl elements is
given by elements $e_{ij}$, with $i,j\in[0;n+m]$; for
an explicit expression of these elements $e_{ij}$ in terms
of the standard Chevalley generators, see~\cite{PalevTolstoy}.
Finding a set of Cartan-Weyl elements, and their (super)commutation
relations, is necessary for the construction of a 
Poincar\'e-Birkhoff-Witt basis of $U_q[gl(n+1|m]$. 
The elements $e_{ij}$ are the $q$-analogues of the defining
basis of $gl(n+1|m)$; their grading is given by $\deg(e_{ij})=
\theta_{ij}$. We shall refer to $e_{ij}$ as a positive root vector
(resp.\ negative root vector) if $i<j$ (resp.\ $i>j$).
For the formulation of the Poincar\'e-Birkhoff-Witt theorem, it is
necessary to fix a total order for the set of elements $e_{ij}$.
Among the positive root vectors, this order is given by
\begin{equation}
e_{ij}<e_{kl},
\hbox{ if } i<k \hbox{ or } i=k \hbox{ and } j<l; 
\label{orderprv}
\end{equation}
for the negative root vectors $e_{ij}$ one takes the same 
rule~(\ref{orderprv}), and total order is fixed by choosing 
\[
\hbox{positive root vectors}<\hbox{negative root vectors}< e_{ii}.
\]
The order among the elements $e_{ii}$ 
is of no real importance since they commute. 
A complete set of relations between the Cartan-Weyl elements
$e_{ij}$ is given by~(\ref{eii-ejj})-(\ref{prv-eij-nrv-ekl})
(see eqs.~(3.10)-(3.15) of~\cite{PalevTolstoy})~:
% ? where should we say something about q ?
\begin{eqnarray}
&& [e_{ii},e_{jj}]=0 ; \label{eii-ejj}\\
&& [e_{ii},e_{jk}]=(\delta_{ij}-\delta_{ik}) e_{jk}; 
\label{eii-ejk}\\
&& (e_{ij})^2=0 \hbox{ if } \t_{ij}=1; \label{eij^2}
\end{eqnarray}
for two positive root vectors $e_{ij}<e_{kl}$~:
\begin{equation}
\lb e_{ij},e_{kl}\rb_{q^{(-1)^{\t_j}\delta_{jl}-(-1)^{\t_j}\delta_{jk}
+(-1)^{\t_i}\delta_{ik}}}
=\delta_{jk}e_{il}+(q-q^{-1})(-1)^{\t_k}\t (l>j>k>i) e_{kj}e_{il};
\label{prv-eij-ekl}
\end{equation}
for two negative root vectors $e_{ij}>e_{kl}$~:
\begin{equation}
\lb e_{ij},e_{kl}\rb_{q^{-(-1)^{\t_j}\delta_{jl}+(-1)^{\t_j}\delta_{jk}
-(-1)^{\t_i}\delta_{ik}}}
=\delta_{jk}e_{il}-(q-q^{-1})(-1)^{\t_k}\t (i>k>j>l) e_{kj}e_{il}; 
\label{nrv-eij-ekl}
\end{equation}
and finally for a positive root vector $e_{ij}$ and a negative
root vector $e_{kl}$~:
\begin{eqnarray}
&&\lb e_{ij},e_{kl}\rb=
\frac{\delta_{il}\delta_{jk}}{{q-q^{-1}}} 
\left(q^{e_{ii}-(-1)^{\t_{ij}}e_{jj}} -q^{-e_{ii}+(-1)^{\t_{ij}}e_{jj}}  
\right)\nn\\
&&+{\Big (} (q-q^{-1}) \theta (j>k>i>l)(-1)^{\t_k}
e_{kj}e_{il} -\delta_{il}\theta (j>k)(-1)^{\t_{kl}}e_{kj} 
+\delta_{jk}\theta (i>l) e_{il}{\Big )} \nn\\
&&\times q^{(-1)^{\t_k}e_{kk}-(-1)^{\t_i}e_{ii}}  
+q^{(-1)^{\t_l}e_{ll}-(-1)^{\t_j}e_{jj}}{\Big (} -(q-q^{-1})
\t (k>j>l>i) (-1)^{\t_j} e_{il}e_{kj} \nn\\
&&-\delta_{il}\theta (k>j)(-1)^{\t_{ij}}e_{kj} 
+\delta_{jk}\theta (l>i) e_{il}{\Big )} .
\label{prv-eij-nrv-ekl}
\end{eqnarray}
Herein,
\begin{equation}
\theta (i_1>i_2>\ldots>i_r)=
\left\{ \begin{array}{ll}
 1, & \hbox{if } i_1>i_2>\ldots>i_r, \\
 0, & \hbox{otherwise.}
\end{array} \right.
\label{def-theta}
\end{equation}

The difference between $U_q[sl(n+1|m)]$ and $U_q[gl(n+1|m)]$ lies
in the elements of the Cartan subalgebra. For $U_q[gl(n+1|m)]$ 
the Cartan subalgebra is generated by $e_{ii}$ ($i\in[0;n+m]$).
For $U_q[sl(n+1|m)]$ the Cartan subalgebra is generated by 
the elements $H_i$, with
\begin{equation}
H_i=e_{00}-(-1)^{\t_i}e_{ii}, \quad i\in[1;n+m].
\end{equation}
Sometimes, it will be useful to work with the elements
$L_i$ and $\bar L_i$, where
\begin{equation}
L_i=q^{H_i},\qquad \bar L_i=q^{-H_i}, \qquad i\in[1;n+m].
\label{LiHi}
\end{equation}
The Cartan-Weyl elements of $U_q[sl(n+1|m)]$ are now given by
$\{ H_i; i\in[1;n+m]\} \cup
\{ e_{ij}; i\ne j \in [0;n+m]\}$. 
The complete set of supercommutation relations between these
Cartan-Weyl elements is given by 
\begin{eqnarray}
&& [H_i,H_j]=0 ; \label{Hi-Hj}\\
&& [H_i,e_{jk}]=(\delta_{0j}-\delta_{0k}-(-1)^{\t_i}
(\delta_{ij}-\delta_{ik})) e_{jk}; \label{Hi-ejk}
\end{eqnarray}
(\ref{eij^2}), (\ref{prv-eij-ekl}), (\ref{nrv-eij-ekl}) and finally
the relation between a positive root vector $e_{ij}$ and a negative
root vector $e_{kl}$~:
\begin{eqnarray}
&&\lb e_{ij},e_{kl}\rb=
\frac{\delta_{il}\delta_{jk}}{{q-q^{-1}}} 
\left( L_j^{(-1)^{\t_i}}{\bar L}_i^{(-1)^{\t_i}}-
{\bar L}_j^{(-1)^{\t_i}} L_i^{(-1)^{\t_i}}
\right)\label{final-rel}\\
&&+ {\Big (} (q-q^{-1}) \theta (j>k>i>l)(-1)^{\t_k}
e_{kj}e_{il}  -\delta_{il}\theta (j>k)(-1)^{\t_{kl}}e_{kj}+
\delta_{jk}\theta (i>l) e_{il}{\Big )} L_i{\bar L}_k \nn \\
&& +L_j{\bar L}_l{\Big (} -(q-q^{-1})
\t (k>j>l>i) (-1)^{\t_j}
e_{il}e_{kj} -\delta_{il}\theta (k>j)(-1)^{\t_{ij}}e_{kj}+
\delta_{jk}\theta (l>i) e_{il}{\Big )} . \nn
\end{eqnarray}

The Jacobson generators of $U_q[sl(n+1|m)]$ are now defined as the
Cartan elements $H_i$ ($i\in[1;n+m]$) together with
the elements
\begin{equation}
a_i^-=e_{0i},\qquad a_i^+=e_{i0},\qquad i\in[1;n+m].
\label{JGs}
\end{equation}

{}From~(\ref{final-rel}) it is easy to deduce that
\begin{equation}
\lb a_i^-, a_j^+ \rb = -(-1)^{\t_i} L_i e_{ji}, \quad (i<j);
\qquad
\lb a_i^-, a_j^+ \rb = -(-1)^{\t_j}e_{ji} \bar L_j, \quad (i>j).
\label{23}
\end{equation}
However, these relations are not complete in order to reshuffle all
Cartan-Weyl elements in an arbitrary expression in the right order. 
For this purpose, we have the following result~:
 
\begin{theo}
A set of Cartan-Weyl elements of $U_q[sl(n+1|m)]$ is given by
$H_i$, $a_i^\pm$, $\lb a_i^+,a_j^- \rb$ ($i\ne j\in[1;n+m]$).
A complete set of supercommutation relations between these 
elements is given by~:
\begin{eqnarray}
&& [H_i,H_j]=0; \qquad
[ H_i,a_{j}^\pm ] =\mp(1+(-1)^{\t_i}\delta_{ij}) a_{j}^\pm ; 
\label{LL} \\
&& 
\lb a_{i}^-, a_i^+\rb={L_i-{\bar{L}_i}\over{q-q^{-1}}}; \label{a-a+}\\
&& \lb a_{i}^\eta, a_j^\eta \rb_q=0\quad (i<j);\qquad
(a^{\pm}_i)^2=0 \quad (i\in[n+1;n+m]);
\label{aa} \\
&& \lb \lb a_i^\eta,a_j^{-\eta}\rb,a_k^\eta\rb_{q^{\xi(1+(-1)^{\t_i}\delta_{ik}})}= \eta^{\t_j}\delta_{jk}L_k^{-\xi\eta}a_i^\eta + 
(-1)^{\t_k}\epsilon(j,k,i)
(q-\q)\lb a_k^\eta,a_j^{-\eta}\rb a_i^\eta \nn\\
&&\qquad = \eta^{\t_j}\delta_{jk}L_k^{-\xi\eta}a_i^\eta + 
(-1)^{\t_k\t_j}\epsilon(j,k,i)q^\xi
(q-\q)a_i^\eta\lb a_k^\eta,a_j^{-\eta}\rb , \label{aaa} \\
&& \hbox{where } (j-i)\xi>0, \;\; \xi,\;\eta =\pm  \nn\\
&& \hbox{and } \epsilon(j,k,i)=
\left\{ \begin{array}{ll}
 1, & \hbox{if }j>k>i;\\
 -1, & \hbox{if }j<k<i;\\
 0, & \hbox{otherwise,}
\end{array}\right. \nn
\end{eqnarray}
and we have used the notation $\bar q = q^{-1}$.
\end{theo}

\noindent {\bf Proof.} 
The first part of the statement is obvious. Relation (\ref{LL})
follows from (\ref{Hi-Hj}) and (\ref{Hi-ejk}); (\ref{a-a+}) 
follows from (\ref{final-rel}) with $l=i<j=k$;
the first relation in (\ref{aa}) follows from (\ref{prv-eij-ekl})
with $i=k<j<l$ and from (\ref{nrv-eij-ekl}) with $l=j<k<i$,
whereas the second relation in (\ref{aa}) comes from (\ref{eij^2}).
Finally, it remains to prove (\ref{aaa}). There are four similar
cases to consider, according to $\eta=\pm$ and $\xi=\pm$. 
For $\eta=-$ and $\xi=+$, we use the first relation
in~(\ref{23}) and find~:
\begin{eqnarray*}
\lb \lb a_i^-,a_j^+\rb,a_k^-\rb_{q^{1+(-1)^{\t_i}\delta_{ik}}}&=&
 -(-1)^{\t_i} \lb L_i e_{ji}, a_k^- \rb_{q^{1+(-1)^{\t_i}\delta_{ik}}}  
= -(-1)^{\t_i} L_i \lb e_{ji}, a_k^- \rb  \\
&=& -(-1)^{\t_i} L_i \lb e_{ji}, e_{0k} \rb = (-1)^{\t_i+\t_{ij}\t_k}
L_i \lb e_{0k},e_{ji} \rb.
\end{eqnarray*}
Herein, we have used the last equation of (\ref{LL}) to change
the order of $L_i$ and $a_k^-$. 
For the last supercommutator, we use (\ref{final-rel})~:
\[
\lb \lb a_i^-,a_j^+\rb,a_k^-\rb_{q^{1+(-1)^{\t_i}\delta_{ik}}} =
(-1)^{\t_i+\t_{ij}\t_k} L_i L_k \bar L_i
\left( -(q-\q)\t(j>k>i)(-1)^{\t_k}e_{0i}e_{jk}+\delta_{kj}e_{0i}\right).
\]
Using trivial properties of the $\t_i$-symbols, the second term 
in the rhs of this expression becomes 
$(-1)^{\t_j} \delta_{jk} L_k a_i^-$; for the
first term we use similar properties and replace according
to (\ref{23}) $e_{jk}$ by $-(-1)^{\t_k}\bar L_k \lb a_k^-,a_j^+\rb$, so
there comes
\[
(-1)^{\t_k}q(q-\q)\t(j>k>i) a_i^-  \lb a_k^-,a_j^+\rb.
\]
This coincides with the second expression in (\ref{aaa}).
Exchanging indices $i$ and $k$, and using the relation just obtained,
one shows that also the first expression in (\ref{aaa}) is
valid. 

For the remaining choices of $\eta$ and $\xi$, the proof is similar.
\mybox

Finally, we wish to remark that in order to construct $U_q[sl(n+1|m)]$
by means of the JGs subject to a set of relations, not all relations of
Theorem~1 are needed.
Such a minimal set of relations was determined in~\cite{PS99}. 

\section{Fock representations}
\setcounter{equation}{0}

In this section we shall construct so-called Fock representations
of $U_q[sl(n+1|m)]$. The representations
considered here are diagonal with respect to the Cartan elements
$H_i$. So it will be convenient to fix $q$ (or $h$) as a
complex number in this and in the following section.

The Fock representations, or modules, can
be defined by means of an induced module construction. 
First observe that $G=U_q[sl(n+1|m)]$, with Cartan-Weyl elements
$H_i$, $a_i^\pm$ and $\lb a_i^+,a_j^- \rb$ ($i\ne j\in[1;n+m]$),
has a subalgebra $H=U_q[gl(n|m)]$ with Cartan-Weyl elements
$H_i$ and $\lb a_i^+,a_j^- \rb$ ($i\ne j\in[1;n+m]$).
A trivial one-dimensional $H$ module is defined as follows~:
\begin{eqnarray}
&& \lb a_i^-,a_j^+\rb|0\ra =0, \qquad (i\ne j\in[1;n+m]) \label{3-1}\\
&& H_i|0\ra =p|0\ra ,
\end{eqnarray}
where $p$ is any complex number. 
Let $P$ be the (associative) subalgebra of $G=U_q[sl(n+1|m)]$
generated by the elements of $H$ and $\{ a^-_i ; i\in[1;n+m]\}$.
The one-dimensional module $\C |0\ra $ can be extended to a
one-dimensional module of $P$ by requiring~:
\begin{equation}
a_i^-|0\ra =0, \qquad i\in[1;n+m]. \label{3-3}
\end{equation}
Now the $G$ module $\bar W_p$ is defined as
\[
\bar W_p= \hbox{Ind}_P^G \;\C |0\ra .
\]
By construction, this means that $\bar W_p$ is freely generated
by the generators $a_i^+$ ($i\in[1;n+m]$) acting on $|0\ra $.
In other words, a basis for $\bar W_p$ is given by
\begin{eqnarray}
&&|p;r_1,r_2,\l ,r_{n+m}\ra  \equiv 
(a_1^+)^{r_1}(a_2^+)^{r_2}\l (a_n^+)^{r_n}(a_{n+1}^+)^{r_{n+1}}
(a_{n+2}^+)^{r_{n+2}} \l (a_{n+m}^+)^{r_{n+m}}|0\ra  
\nn \\
&&\hbox{where }r_i\in \Z_+ \hbox{ for }i\in[1;n]
\hbox{ and }r_{i}\in \{0,1\}\hbox{ for }i\in[n+1;m]. \label{def-vector}
\end{eqnarray} 
So $\bar W_p$ is an infinite-dimensional $G$ module. 
The main part of this section is devoted to the computation
of the action of the JGs on the basis vectors~(\ref{def-vector})
of $\bar W_p$. This, of course, completely determines the action
of $U_q[sl(n+1|m)]$ on $\bar W_p$.  

\begin{theo}
The transformation of the basis~(\ref{def-vector}) of $\bar W_p$
under the action of the JGs reads~:
\begin{eqnarray}
&&H_i|p;r_1,r_2,\ldots,r_{n+m}\ra =\left(p-(-1)^{\t_i}r_i-\sum_{j=1}^{n+m}r_j
\right) |p;r_1,r_2,\ldots,r_{n+m}\ra , \label{actionH}\\
&&a_{i}^-|p;r_1,r_2,\ldots, r_{n+m}\ra = 
(-1)^{\t_1r_{1}+\t_2r_2 +\l +\t_{i-1}r_{i-1}}
q^{r_1+\l +r_{i-1}}[r_i] [p-\sum_{j=1}^{n+m}r_j+1]  \nn\\
&&\qquad \times |p;r_1,r_2,\ldots, r_{i-1},r_i-1,r_{i+1},\l 
, r_{n+m}\ra , \label{actiona-}\\[3mm]
&& a_{i}^+|p;r_1,r_2,\ldots, r_{n+m}\ra = 
(-1)^{\t_1r_{1}+\t_2r_2 +\l +\t_{i-1}r_{i-1}}
 \q^{r_1+\l+r_{i-1}} (1-\t_i r_i) \nn\\
&&\qquad\times |p;r_1,r_2,\ldots, r_{i-1},r_i+1,r_{i+1},\l , r_{n+m}\ra,
\label{actiona+}
\end{eqnarray}
where $i\in[1;n+m]$.
\label{theoaction}
\end{theo}

\noindent {\bf Proof.}
Equation~(\ref{actionH}) is an immediate consequence of
$[H_i,a_j^+] = - (1+(-1)^{\t_i}\delta_{ij})a_j^+$,
which is one of the last relations in~(\ref{LL}).
Also the action of $a_i^+$ on the basis vectors is easy~:
(\ref{actiona+}) follows directly from~(\ref{aa}).
The hard work lies in the proof of~(\ref{actiona-}).
For this purpose, we shall use a number of technical lemmas.

\begin{lemm}
The following relations hold~:
\begin{eqnarray}
\bullet&&
\lb A,B_1B_2\ldots B_{i-1}B_{i}B_{i+1}\ldots B_j\rb_{q^{b_1+b_2+
\ldots+b_j}} \nn\\
&& =\sum_{i=1}^j
{q^{b_1+b_2+\ldots+b_{i-1}}}(-1)^{\alpha(\beta_1+\ldots +\beta_{i-1})}
B_1B_2\ldots B_{i-1}\lb A,B_{i}\rb_{q^{b_i}}B_{i+1}\ldots B_j ,
\nn \\
&& \qquad\hbox{where }\alpha=\deg(A)\hbox{ and }\beta_i=\deg(B_i); 
\label{45} \\
\bullet&& \lb a_i^-,(a_j^+)^r\rb=
\left\{\begin{array}{ll}
\ds \frac{\q^{2r}-1}{\q^2-1}(a_j^+)^{r-1} \lb a_i^-,a_j^+\rb &
\hbox{ when } i<j, \\[3mm]
\ds \frac{q^{2r}-1}{q^2-1}(a_j^+)^{r-1}\lb a_i^-,a_j^+\rb &
\hbox{ when } i>j; 
\end{array}\right. \label{46}\\
\bullet&& \lb a_i^-,(a_i^+)^r \rb=\frac{(a_i^+)^{r-1}}{q-\q}
\left( \frac{\q^{2r}-1}{\q^2-1} L_i-
\frac{q^{2r}-1}{q^2-1} \L_i \right); \label{47} \\
\bullet&&
\lb \lb a_i^-,a_j^+\rb,(a_i^+)^r \rb_{q^r}=
-(-1)^{\t_j}
\frac{\q^{2r}-1}{\q^2-1} \L_i a_j^+ (a_i^+)^{r-1}, \quad i>j,
\label{48} \\
\bullet&& \lb \lb a_i^-,a_j^+\rb,(a_k^+)^{r}\rb_{q^{r}}=
(-1)^{\t_j}(q^{2{r}}-1)a_j^+(a_k^+)^{{r}-1}\lb a_i^-,a_k^+\rb,
\quad i>k>j. \label{49}
\end{eqnarray}
\end{lemm} 

\noindent {\bf Proof.} 
Equation~(\ref{45}) follows by direct calculation. 
We need to prove equation~(\ref{46}) only when $r>1$, i.e.\ only
when $\t_j=\bar 0$. Then one writes, using~(\ref{45}),
\[
\lb a_i^-,(a_j^+)^r\rb= \lb a_i^-,(a_j^+)^{r-1} a_j^+ \rb
= \lb a_i^-,(a_j^+)^{r-1}\rb a_j^+ + (a_j^+)^{r-1} \lb a_i^-, a_j^+ \rb .
\]
Now the result follows using induction on $r$ and using the 
triple relation~(\ref{aaa}) with $k=j$ and $\eta=-$.
The proof of~(\ref{47}) is similar, using~(\ref{45}), induction
on $r$, and~(\ref{a-a+}).
Also the proof of~(\ref{48}) goes along the same line~: 
first one writes $(a_i^+)^r$ as $a_i^+(a_i^+)^{r-1}$
(for $r>1$); using~(\ref{45}) this yields two terms~: on
the first term one applies~(\ref{aaa}), and on the second term
one applies~(\ref{48}) by induction; then the result follows.
The proof of~(\ref{49}) is essentially the same.
\mybox

\begin{lemm}
For $i>1$ the following relation holds~:
\begin{eqnarray}
&& \lb a_i^-,a_1^+\rb(a_2^+)^{r_2}\ldots (a_{n+m}^+)^{r_{n+m}}|0\ra \nn\\
&& =   - (-1)^{\t_1+\t_2r_2+\t_3r_3+\l +\t_{i-1}r_{i-1}} 
q^{2r_2+\ldots+2r_{i-1}+r_i+\ldots+r_{n+m}-p} [r_i] \nn\\
&& \times a_1^+ (a_2^+)^{r_2}\ldots (a_{i-1}^+)^{r_{i-1}}
(a_i^+)^{r_i-1}(a_{i+1}^+)^{r_{i+1}}\ldots(a_{n+m}^+)^{r_{n+m}}|0\ra.
\label{51a}
\end{eqnarray}
\end{lemm}

\noindent
{\bf Proof.}
Consider first $i=2$. Using~(\ref{45}), one finds
\begin{eqnarray}
&& \lb a_2^-,a_1^+\rb(a_2^+)^{r_2}\ldots (a_{n+m}^+)^{r_{n+m}}|0\ra =
\lb \lb a_2^-,a_1^+\rb, (a_2^+)^{r_2}\ldots (a_{n+m}^+)^{r_{n+m}}
\rb_{q^{r_2+\l +r_{n+m}}}|0\ra \nn\\
&&=  \lb \lb a_2^-,a_1^+\rb, (a_2^+)^{r_2}\rb_{q^{r_2}} 
(a_3^+)^{r_3}\ldots (a_{n+m}^+)^{r_{n+m}}|0\ra  \nn\\
&& +(-1)^{(\t_1+\t_2)r_2\t_2}q^{r_2}(a_2^+)^{r_2}
\lb \lb a_2^-,a_1^+\rb, (a_3^+)^{r_3}\ldots (a_{n+m}^+)^{r_{n+m}}
\rb_{q^{r_3+\l +r_{n+m}}}|0\ra .
\label{tmp}
\end{eqnarray}
{}From~(\ref{45}) and~(\ref{aaa}) it follows that the second term in the
rhs of~(\ref{tmp}) is zero. For the first term, apply~(\ref{48}) and
use the action of $\bar L_2$ as given by~(\ref{actionH}) and~(\ref{LiHi}).
Then the result follows.\\
Next we shall use induction on $i$ to prove~(\ref{51a}) in general.
So suppose~(\ref{51a}) holds for all $j=2,3,\ldots,i-1$, i.e.
\begin{eqnarray}
&& \lb a_j^-,a_1^+\rb(a_2^+)^{r_2}\ldots (a_{n+m}^+)^{r_{n+m}}|0\ra \nn\\
&&= - (-1)^{\t_1+\t_2r_2+\t_3r_3+\l +\t_{j-1}r_{j-1}} 
 q^{2r_2+\ldots+2r_{j-1}+r_j+\ldots+r_{n+m}-p} [r_j]  \nn\\
&& \times a_1^+ (a_2^+)^{r_2}\ldots (a_{j-1}^+)^{r_{j-1}}
  (a_j^+)^{r_j-1}(a_{j+1}^+)^{r_{j+1}}\ldots(a_{n+m}^+)^{r_{n+m}}|0\ra .
\label{56}
\end{eqnarray}
Making a shift of indices in~(\ref{56}) (thereby putting the last
$r_k$-values equal to zero), leads to the following equivalent
equation~:
\begin{eqnarray}
&& \lb a_{i}^-,a_j^+\rb(a_{j+1}^+)^{r_{j+1}}\ldots 
(a_{n+m}^+)^{r_{n+m}}|0\ra \nn\\
&&= - (-1)^{\t_j+\t_{j+1}r_{j+1}+\t_{j+2}r_{j+2}+\l +\t_{i-1}r_{i-1}} 
 q^{2r_{j+1}+\ldots+2r_{i-1}+r_{i}+ \ldots+r_{n+m}-p} [r_{i}] \nn\\
&& \times a_j^+ (a_{j+1}^+)^{r_{j+1}}\ldots (a_{i-1}^+)^{r_{i-1}}
  (a_{i}^+)^{r_{i}-1}(a_{i+1}^+)^{r_{i+1}}
\ldots(a_{n+m}^+)^{r_{n+m}}|0\ra , \quad (j<i).
\label{59}
\end{eqnarray}
Now consider the lhs of~(\ref{51a}) and apply~(\ref{45})~:
\begin{eqnarray}
&& \lb a_i^-,a_1^+\rb(a_2^+)^{r_2}\ldots (a_{n+m}^+)^{r_{n+m}}|0\ra =
\lb \lb a_i^-,a_1^+\rb, (a_2^+)^{r_2}\ldots (a_{n+m}^+)^{r_{n+m}}
\rb_{q^{r_2+\l +r_{n+m}}}|0\ra \nn\\
&& = \sum_{k=2}^{n+m}q^{r_2+\l +r_{k-1}}
(-1)^{(\t_i+\t_1)(\t_2r_2+\t_3r_3+\l +\t_{k-1}r_{k-1})} & \cr
&& \times (a_2^+)^{r_2}\ldots (a_{k-1}^+)^{r_{k-1}}
\lb \lb  a_i^-, a_1^+\rb, 
  (a_k^+)^{r_k}\rb_{q^{r_k}}(a_{k+1}^+)^{r_{k+1}}
\ldots(a_{n+m}^+)^{r_{n+m}}|0\ra.
\label{tmp2}
\end{eqnarray}
In this last sum, all terms with $k>i$ are easily seen to vanish.
For the terms with $k<i$, we apply~(\ref{49}), and for the term
with $k=i$, we apply~(\ref{48}). Then there comes~:
\begin{eqnarray}
&& \sum_{k=2}^{i-1} (-1)^{\t_1+\t_2r_2+\t_3r_3+\l 
+\t_{k-1}r_{k-1} }  (q^{2r_k}-1) \nn\\
&& \qquad \times a_1^+(a_2^+)^{r_2}\ldots (a_{k-1}^+)^{r_{k-1}}
(a_k^+)^{r_k-1} \lb  a_i^-, a_k^+\rb 
(a_{k+1}^+)^{r_{k+1}} \ldots(a_{n+m}^+)^{r_{n+m}}|0\ra \nn\\
&& -(-1)^{\t_1+\t_2r_2+\t_3r_3+\l +\t_{i-1}r_{i-1} } 
q^{-p+\sum_{l=i}^{n+m}r_l}[r_i] \nn\\
&& \qquad \times a_1^+(a_2^+)^{r_2}\ldots (a_{i-1}^+)^{r_{i-1}} 
(a_i^+)^{r_i-1}(a_{i+1}^+)^{r_{i+1}}
\ldots(a_{n+m}^+)^{r_{n+m}}|0\ra.  
\label{tmp3}
\end{eqnarray}
For the terms in~(\ref{tmp3}) with $k<i$, we can apply~(\ref{59}).
Then it is a matter of appropriately summing all contributions,
which leads finally to the rhs of~(\ref{51a}). \mybox

\noindent
{\bf Proof of Theorem~\ref{theoaction}.}
There remains to prove equation~(\ref{actiona-}). 
First, assume that $i=1$ in~(\ref{actiona-});
then we have according to~(\ref{45})
\begin{eqnarray}
&& a_{1}^-|p;r_1,r_2,\ldots, \l , r_{n+m}\ra = 
\lb a_1^-, (a_1^+)^{r_1}(a_2^+)^{r_2}\l 
 (a_{n+m}^+)^{r_{n+m}}\rb|0\ra \nn\\
&=& \lb a_1^-, (a_1^+)^{r_1}\rb (a_2^+)^{r_2}\l 
 (a_{n+m}^+)^{r_{n+m}}|0\ra + 
\sum_{j=2}^{n+m} (-1)^{\t _1(\t_1r_1+\t_2r_2+\l +\t_{j-1}r_{j-1})} \nn\\
&& \times(a_1^+)^{r_1}\l (a_{j-1}^+)^{r_{j-1}} \lb a_1^-, (a_j^+)^{r_j}\rb
(a_{j+1}^+)^{r_{j+1}}\l 
 (a_{n+m}^+)^{r_{n+m}}|0\ra . \label{tmp4}
\end{eqnarray}
The terms with $j\geq 2$ in the rhs of~(\ref{tmp4}) are found to
be zero using~(\ref{46}) and~(\ref{aaa}). So only the first term
in the rhs of~(\ref{tmp4}) gives a contribution; using~(\ref{47})
this is
\[
[r_1] [p-\sum_{j=1}^{n+m}r_j+1] |p;r_1-1,r_2,\ldots, r_{n+m}\ra,
\]
so the case $i=1$ is proved. Now we use again induction on $i$.
So the following equation holds for $j<i$~:
\begin{eqnarray}
&& a_{j}^- (a_1^+)^{r_1}\l (a_{n+m}^+)^{r_{n+m}}|0\ra= 
(-1)^{\t_1r_{1}+\t_2r_2 +\l +\t_{j-1}r_{j-1} }
q^{r_1+\l +r_{j-1}} \nn\\
&&\quad\times [r_j] [p-\sum_{l=1}^{n+m}r_l+1]  
|p;r_1,r_2,\ldots, r_{j-1},r_j-1,r_{j+1},\l , r_{n+m}\ra .
\label{52}
\end{eqnarray}
In this equation, put $r_{n+m}=0$ and raise all indices by~1.
Then the following (equivalent) equation holds~:
\begin{eqnarray}
&& a_{i}^-(a_2^+)^{r_2}\l (a_{n+m}^+)^{r_{n+m}}|0\ra = 
(-1)^{\t_2r_{2}+\t_3r_3 +\l +\t_{i-1}r_{i-1} }
q^{r_2+\l +r_{i-1}} \label{54}\\
&&\times [r_{i}] [p-\sum_{l=2}^{n+m}r_l+1]  
(a_2^+)^{r_2}\l (a_{i-1}^+)^{r_{i-1}}
(a_{i}^+)^{r_{i}-1}(a_{i+1}^+)^{r_{i+1}}\l (a_{n+m}^+)^{r_{n+m}}
|0\ra . \nn
\end{eqnarray}
Now consider
\begin{eqnarray}
&& a_{i}^-|p;r_1,r_2,\l ,  r_{n+m}\ra = 
\lb  a_i^-, (a_1^+)^{r_1}\l (a_{n+m}^+)^{r_{n+m}}\rb|0\ra \nn\\
&& = \lb  a_i^-, (a_1^+)^{r_1}\rb (a_2^+)^{r_2}\l  
(a_{n+m}^+)^{r_{n+m}} \rb |0\ra + (-1)^{\t _i\t_1r_{1}} 
(a_1^+)^{r_1}\lb  a_i^-,  (a_2^+)^{r_2}\l  
(a_{n+m}^+)^{r_{n+m}}\rb |0\ra \nn\\
&& =  {q^{2r_1-1}\over{q^2-1}}(a_1^+)^{r_1-1}\lb  a_i^-, a_1^+\rb 
(a_2^+)^{r_2}\l (a_{n+m}^+)^{r_{n+m}}\rb|0\ra  \nn\\
&&\quad + (-1)^{\t _i\t_1r_{1}} 
(a_1^+)^{r_1}\lb  a_i^-,  (a_2^+)^{r_2}\l  
(a_{n+m}^+)^{r_{n+m}}\rb |0\ra.
\label{tmp5}
\end{eqnarray}
This was obtained by applying~(\ref{46}) on the first term.
The rhs of~(\ref{tmp5}) can now be determined as follows~:
for the first term we use~(\ref{51a}),
and for the second term, we use~(\ref{54}) in which both
sides have been multiplied (on the left) by $(a_1^+)^{r_1}$.
Adding both contributions leads to the desired result.
\mybox
 
The action of the elements $H_i$ and $a_i^\pm$ ($i\in[1;n+m]$)
on the basis
vectors of $\bar W_p$, determined in Theorem~\ref{theoaction},
clearly imply that $\bar W_p$ has an invariant submodule
when $p$ is a nonnegative integer. From now on we shall assume
that $p\in\Z_+$. Then we have

\begin{coro}
The $U_q[sl(n+1|m)]$ module $\bar W_p$ has an invariant
submodule $V_p$ with basis vectors
\[
|p;r_1,r_2,\ldots,r_{n+m}\ra, \hbox{ with }
\sum_{i=1}^{n+m} r_i >p.
\]
The quotient module $W_p = \bar W_p / V_p$ is an irreducible
representation for $U_q[sl(n+1|m)]$. The basis vectors of
$W_p$ are given by (the representatives of)
\begin{equation}
|p;r_1,r_2,\ldots,r_{n+m}\ra, \hbox{ with }
\sum_{i=1}^{n+m} r_i \leq p.
\label{basisofWp}
\end{equation}
\end{coro}

These finite-dimensional irreducible $U_q[sl(n+1|m)]$ modules
$W_p$ are referred to as the Fock modules or Fock representations
of $U_q[sl(n+1|m)]$. Also in the Fock modules, the action of the
elements $H_i$ and $a_i^\pm$ ($i\in[1;n+m]$) on the basis 
vectors~(\ref{basisofWp}) is essentially given by the equations
of Theorem~\ref{theoaction}.

One can verify that the irreducible Fock representations
$W_p$ are so-called atypical representations of $U_q[sl(n+1|m)]$.
Atypicality is usually defined for highest weight representations
of simple Lie superalgebras~\cite{Kac2}, 
but it can be extended to highest weight
representations of the corresponding Hopf superalgebras~\cite{Zhang}.
In the standard basis, the Dynkin labels of $W_p$ (or of
its highest weight) are given by $(p,0,\ldots,0)$. 
This means that in general the representation $W_p$ is
multiply atypical~\cite{Kac2,VHKT}. 
More precisely, if $n\geq m$, then $W_p$ is
$m$-fold atypical; if $n<m$, then $W_p$ is $(n+1)$-fold atypical for
$p<m-n$ and $n$-fold atypical for $p\geq m-n$. 
Observe that in this way we have obtained the action of 
a set of generators of $U_q[sl(n+1|m)]$ on a class
of atypical irreducible representations, i.e.\ the Fock modules.
In general, an explicit basis for atypical representations is not known,
not even in the case of $sl(n+1|m)$. 
For typical representations of $U_q[sl(n+1|m)]$, it is easier 
to construct a basis. 
For a subclass of these, the so-called essentially typical
representations, a (Gelfand-Zetlin)
basis has been constructed together with the action of the Chevalley
generators~\cite{PSV94}.

\section{Unitary Fock representations}
\setcounter{equation}{0}

In this section we select a class of Fock modules important
for physical applications. These are the ones for which the
standard Fock metric is positive definite, and for which the
representatives of $a_i^\pm$ and $H_i$ ($i\in[1;n+m]$) satisfy
the Hermiticity conditions~:
\begin{equation}
(a_i^+)^\dagger = a_i^-,\qquad
(a_i^-)^\dagger = a_i^+,\qquad
(H_i)^\dagger = H_i.
\label{dagger}
\end{equation}
In quantum mechanics, including its generalization
to the noncommutative case (see, for instance~\cite{Kempf,PN97}),
(\ref{dagger}) follows from the relations 
$a_k^\pm =\hbox{const} (x_k \mp i p_k)$
and the requirement that the position operators $x_k$
and the momentum operators $p_k$ should be selfadjoint operators.
By definition, representations for which~(\ref{dagger}) holds
are said to be unitary
(with respect to the anti-involution in $U_q[sl(n+1|m)]$
defined by~(\ref{dagger}), and the Fock space scalar product).

For the Fock representation $W_p$, we can define a Hermitian
form $(\,,\,)$ by requiring
\begin{equation}
( |0\ra\,,\, |0\ra) = \la 0 | 0 \ra =1,
\label{00}
\end{equation}
and by postulating that the Hermiticity conditions~(\ref{dagger})
should be satisfied, i.e.
\begin{equation}
( a_i^\pm v, w ) = (v, a_i^\mp w),\qquad
\forall v, w \in W_p.
\label{Hermcond}
\end{equation}
It is now easy to determine that any two vectors
$|p;r_1,r_2,\ldots,r_{n+m}\ra$ and $|p;r_1',r_2',\ldots,r_{n+m}'\ra$
with $(r_1,r_2,\ldots,r_{n+m}) \ne (r_1',r_2',\ldots,r_{n+m}')$
are orthogonal. Furthermore, one can compute~:
\begin{equation}
\left(|p;r_1,r_2,\ldots,r_{n+m}\ra,|p;r_1,r_2,\ldots,r_{n+m}\ra\right)=
 {[p]!\over [p-R]!}\prod_{i=1}^{n+m} [r_i]!=
 {[p]!\over [p-R]!}\prod_{i=1}^n [r_i]!,
\label{squarednorm}
\end{equation}
where $R=r_1+r_2+\ldots+r_{n+m}$.
Clearly, it holds for $R=0$; then use induction on~$R$ together
with~(\ref{actiona-})-(\ref{actiona+}).

Assume now that $1\leq i < j \leq n+m$. According
to~(\ref{squarednorm}) we have 
\begin{equation}
( a_i^+a_j^+ |0\ra\,,\, a_i^+a_j^+ |0\ra ) = [p][p-1].
\label{p-1}
\end{equation}
{}From~(\ref{aa}) we have $a_i^+a_j^+=(-1)^{\t_i\t_j}
q a_j^+a_i^+ = (-1)^{\t_i} q a_j^+a_i^+$ (since
$(-1)^{\t_i\t_j}=(-1)^{\t_i}$ for $i<j$); thus we find
\begin{equation}
( a_j^+a_i^+ |0\ra\,,\, a_i^+a_j^+ |0\ra ) =
( (-1)^{\t_i}\bar q a_i^+a_j^+ |0\ra\,,\, a_i^+a_j^+ |0\ra ) =
(-1)^{\t_i} {\bar q}^* [p][p-1],
\label{a1a2-1}
\end{equation}
where ${\bar q}^*$ is the complex conjugate of ${\bar q}=q^{-1}$.
On the other hand, using~(\ref{actiona-}),
\begin{eqnarray}
( a_j^+a_i^+ |0\ra\,,\, a_i^+a_j^+ |0\ra ) &=&
( a_i^+ |0\ra\,,\, a_j^- |p; 0,\ldots,0,1_i,0,\ldots,0,1_j,0,\ldots,0\ra) \nn\\
&=&
( a_i^+ |0\ra\,,\, (-1)^{\t_i} q [p-1] |p;0,\ldots,0,1_i,0,\ldots,0\ra  ) \nn\\
&=& (-1)^{\t_i} q [p-1] ( a_i^+ |0\ra\,,\, a_i^+ |0\ra ) 
= (-1)^{\t_i} q[p][p-1].
\label{a1a2-2}
\end{eqnarray}
Herein, $1_i$ stands for a number~1 at the position~$i$.
When $p\geq 2$, the comparison of~(\ref{a1a2-1})
and~(\ref{a1a2-2}) yields $|q|^2=1$. Hence a necessary condition
for the Fock space to be unitary is that $q$ must be a phase, i.e.
\begin{equation}
q= e^{i\phi}, \qquad (-\pi < \phi < \pi).
\label{phase}
\end{equation}

Let us now further investigate when the Fock module is unitary,
i.e.\ when the Hermitian form $(\,,\,)$ is an inner product.
This means that for every $(r_1,\ldots,r_{n+m})$ 
with $0\leq R\leq p$, the value
in~(\ref{squarednorm}) should be positive. 
In particular, this implies that all the numbers
\[
[p],\; [p-1],\; [p-2],\; \ldots,\; [2],\; [1]
\]
should be positive.
However, since $q=e^{i\phi}$ is a phase, we have
\[
[k] = {q^k-q^{-k} \over q-q^{-1}} = {\sin(k\phi)\over\sin(\phi)}.
\]
So we are left with the following question~: let $p>1$, find the
values of $\phi$ ($-\pi<\phi<\pi$) where all of the following
functions
\[
{\sin(2\phi)\over \sin(\phi)}, {\sin(3\phi)\over \sin(\phi)},
\ldots, {\sin(p\phi)\over \sin(\phi)}
\]
are positive. For each of these functions ${\sin(k\phi)\over \sin(\phi)}$,
the zeros and hence the signs are easy to determine. 
So the common domain where all of these functions are positive
is given by
\[
{-\pi\over p} < \phi < {\pi\over p}.
\]
Thus we have
\begin{theo}
The irreducible Fock module $W_p$ ($p\geq 2$) is 
unitary if and only if $q$ is a 
phase, i.e.\ $q=e^{i\phi}$, with ${-\pi\over p} < \phi < {\pi\over p}$.
\label{theounitary}
\end{theo}

Observe that whether $q$ is a root of unity or not does not have 
any effect on the irreducibility or unitarity of the Fock module
$W_p$, as long as the conditions of Theorem~\ref{theounitary}
are satisfied. Indeed, suppose that $q=e^{i\phi}$ is a root of unity
with $\phi$ a rational multiple of $\pi$ and ${-\pi\over p} < \phi 
< {\pi\over p}$. Then the smallest integer $N$ for which $q^N=-1$
is greater than $p$. As a consequence, the rhs in~(\ref{squarednorm})
is never zero. This implies that there are no singular vectors
among the weight vectors $|p;r_1,\ldots, r_{n+m}\ra$, and
thus irreducibility holds.

Under the conditions of Theorem~\ref{theounitary}, 
we can define an orthonormal basis of $W_p$~:
\begin{equation}
|p;r_1,r_2,\ldots,r_{n+m})=\sqrt{[p-\sum_{l=1}^{n+m} r_l]!
\over {[p]![r_1]!\ldots[r_{n+m}]!}}|p;r_1,r_2,\ldots,r_{n+m}\ra ,
\label{64}
\end{equation}
where $0\leq \sum_{i=1}^{n+m}r_i \leq p$. 
In the new basis~(\ref{64}) the transformation 
formulas~(\ref{actionH})-(\ref{actiona+}) read ($i\in[1;n+m]$)~:
\begin{eqnarray}
&&H_i|p;r_1,r_2,\ldots,r_{n+m}) =
\left(p-(-1)^{\t_i}r_i-\sum_{j=1}^{n+m}r_j
\right) |p;r_1,r_2,\ldots,r_{n+m}) , \label{actionHi}\\
&& a_i^-|p;r_1,\ldots,r_{n+m})=
(-1)^{\t_1r_1+\l +\t_{i-1}r_{i-1} } \label{65}\\
&& \times q^{r_1+\ldots +r_{i-1}}\sqrt{[r_i]
[p-\sum_{l=1}^{n+m} r_l +1]}\;\;|p;r_1,\ldots r_{i-1},r_i-1,r_{i+1},
\ldots,r_{n+m}),
\nn\\
&&  a_i^+|p;r_1,\ldots,r_{n+m})=
(-1)^{\t_1r_1+\l +\t_{i-1}r_{i-1} } \label{66}\\
&& \times \q^{r_1+\ldots +r_{i-1}}(1-\t_ir_i)
\sqrt{[r_i+1][p-\sum_{l=1}^{n+m} r_l]}
\;\;|p;r_1,\ldots r_{i-1},r_i+1,r_{i+1},\ldots,r_{n+m}).
\nn
\end{eqnarray}
{}From~(\ref{23}) it is now easy to determine the action of the remaining 
Cartan-Weyl generators $e_{ji}$ on the basis elements of $W_p$~:
\begin{eqnarray}
&& e_{ji}|p;r_1,\ldots,r_{n+m})= (-1)^{\t_i(r_i+1)+\t_{i+1}r_{i+1}+ 
\l +\t_{j-1}r_{j-1}}
{\bar q}^{r_{i+1}+\ldots +r_{j-1}-2\t_i(1-r_i)} (1-\t_jr_j) \nn\\
&& \times \sqrt{[r_i][r_j+1]}\;\;
|p;r_1,\ldots r_{i-1},r_i-1,r_{i+1},
\ldots,r_{j-1},r_j+1,r_{j+1}, \ldots , r_{n+m}),\quad(i<j),\nn\\ 
&& \label{ejip}\\
&& e_{ji}|p;r_1,\ldots,r_{n+m})=  (-1)^{\t_jr_j+ \l +\t_{i-1}r_{i-1}}
q^{2\t_jr_j+r_{j+1}+\ldots +r_{i-1}} (1-\t_jr_j) \nn\\
&& \times \sqrt{[r_i][r_j+1]}\;\;
|p;r_1,\ldots r_{j-1},r_j+1,r_{j+1},
\ldots,r_{i-1},r_i-1,r_{i+1}, \ldots , r_{n+m}),\quad(i>j). \nn\\
&&\label{eijp}
\end{eqnarray}
In particular, it is possible to extend $W_p$ to a $U_q[gl(n+1|m)]$
module, the actions being given by~(\ref{65}), (\ref{66}), 
(\ref{ejip}), (\ref{eijp}) and
\begin{eqnarray}
&& e_{00}|p;r_1,\ldots,r_{n+m})=
(p-\sum_{l=1}^{n+m}r_l)|p;r_1,\ldots,r_{n+m}), \label{82a} \\
&& e_{ii}|p;r_1,\ldots,r_{n+m})=
r_i|p;r_1,\ldots,r_{n+m}), \;\; i\in [1;n+m]. \label{82b}
\end{eqnarray}

\section{Dyson and Holstein-Primakoff realizations of $U_q[sl(n+1|m)]$}
\setcounter{equation}{0}

Consider $(n+m)$ $\Z_2$-graded indeterminates $c_i^\pm$
($i\in[1;n+m]$) with
\begin{equation}
\deg(c_i^\pm)=\t_i.
\end{equation}
Denote by $W(n|m)$ the free $\C[[h]]$ module (completed
in the $h$-adic topology)
generated by the elements $c_i^\pm$ subject to the relations
\begin{equation}
\lb  c_i^-, c_j^+\rb =\delta_{ij}, \qquad \lb  c_i^+, c_j^+\rb =
\lb  c_i^-, c_j^-\rb =0. \label{cc}
\end{equation}
As usual, let
\begin{equation}
N_i=c_i^+c_i^-, \qquad N=\sum_{j=1}^{n+m}N_j. \label{Ni}
\end{equation}
The algebra $W(n|m)$ of $n$ pairs of Bose and $m$ pairs of Fermi
CAO's has a natural action in the Fock space ${\cal F}(n|m)$,
defined as follows. Let ${\cal F}(n|m)$ be the free $W(n|m)$ module
generated by a vector $|0\ra$ subject to the relations
\[
c_i^- |0\ra =0, \hbox{ for all } i\in[1;n+m].
\]
Then it follows easily that a basis of ${\cal F}(n|m)$ is given by
\begin{equation}
(c_1^+)^{l_1}(c_2^+)^{l_2}\l (c_{n+m}^+)^{l_{n+m}}|0\ra\equiv
|l_1,l_2,\l ,l_{n+m}\ra, \label{basisFock}
\end{equation}
where
\[
l_i\in\Z_+\hbox{ for }i\in[1;n]\hbox{ and }l_i\in\{ 0,1\}
\hbox{ for }i\in[n+1;n+m].
\]

The Dyson~\cite{D} and Holstein-Primakoff~\cite{HP} 
realizations of $U_q[sl(n+1|m)]$
are two different algebra homomorphisms of $U_q[sl(n+1|m)]$ 
into $W(n|m)$~\cite{Palev99}. 
Since $W(n|m)$ has the natural Fock representation
${\cal F}(n|m)$, these realizations will provide representations
of $U_q[sl(n+1|m)]$ in ${\cal F}(n|m)$.

\begin{theo}[Dyson realization]
Let $p$ be any complex number. 
The linear map $\rho : U_q[sl(n+1|m)] \rightarrow W(n|m)$, 
defined on the Jacobson generators by
\begin{eqnarray}
\rho (H_i) &=& p-(-1)^{\t_i}c_i^+c_i^- -\sum_{j=1}^{n+m}c_j^+c_j^-
=p-(-1)^{\t_i}N_i-N, \nn\\
\rho (a_i^-)&=&q^{N_1+\l +N_{i-1}}{[N_i+1]\over{N_i+1}}
[p-N]c_i^-, \nn\\
\rho (a_i^+)&=&\bar{q}^{N_1+\l +N_{i-1}}c_i^+, \label{dyson}
\end{eqnarray}
is a (associative algebra) homomorphism of $U_q[sl(n+1|m)]$ into $W(n|m)$.
\label{propdyson}
\end{theo}

The inspiration of this mapping comes from Theorem~\ref{theoaction}.
The actual proof of Theorem~\ref{propdyson} is straightforward
but tedious~:
one has to verify that all relations in Theorem~1 are satisfied
under the substitution of $H_i$ and $a_i^\pm$ by 
$\rho (H_i)$ and $\rho (a_i^\pm )$. These computations are
lengthy and based upon easy relations such as
\[
f(N_i)c_j^\pm =c_j^\pm f(N_i \pm \delta_{ij}),\quad i,j\in[1;n+m];\qquad
q^{N_i}=1-N_i+qN_i \hbox{ for }i>n,
\]
or simple $q$-identities such as $[x+1][y]-[x][y+1]=[y-x]$.

The Dyson realization of the JGs of $U_q[sl(n+1|m)]$ leads
to an explicit realization of all Cartan-Weyl elements
of $U_q[sl(n+1|m)]$ in terms of the Bose and Fermi CAO's.
Indeed, using~(\ref{23}) and~(\ref{dyson}) one obtains~:
\begin{eqnarray}
\rho(e_{ji})&=&q^{2\t_j(N_j-1)+N_{j+1}+N_{j+2}+\l +N_{i-1}}
{[N_i+1]\over{N_i+1}}c_j^+c_i^-, \quad (j<i) \label{dysoneij}\\
\rho(e_{ji})&=&{\bar {q}}^{2\t_iN_i+N_{i+1}+N_{i+2}+\l +N_{j-1}}
{[N_i+1]\over{N_i+1}}c_j^+c_i^-, \quad (j>i). \label{dysoneji}
\end{eqnarray}
In~(\ref{dysoneij}), the convention is that the summation (in the
power of $q$) is 0 when $j=i-1$ (and similarly for~(\ref{dysoneji})).
Since ${\cal F}(n|m)$ is a $W(n|m)$ module, the Dyson realization
provides a representation of $U_q[sl(n+1|m)]$ into ${\cal F}(n|m)$.
It is easy to see that the action of every $\rho(H_i)$ and
$\rho(a_i^\pm)$ upon $|l_1,\ldots,l_{n+m}\ra$ is the same
as the action of $H_i$ and $a_i^\pm$ in the representation
on $\bar W_p$ given by Theorem~\ref{theoaction}, under
the identification 
\[
|l_1,\ldots,l_{n+m}\ra \equiv |p; l_1,\ldots,l_{n+m}\ra.
\]
Therefore, it follows that the representation $\rho$ of $U_q[sl(n+1|m)]$
into ${\cal F}(n|m)$ (under the Dyson realization) is irreducible
when $p\not\in \Z_+$. When $p\in\Z_+$, the representation $\rho$
is indecomposable. The subspace ${\cal F}_1(n|m)$, spanned
on the vectors
\[
|l_1,\ldots,l_{n+m}\ra \hbox{ with } l_1+\cdots +l_{n+m}>p
\]
is clearly invariant under the action of $U_q[sl(n+1|m)]$.
We denote the (finite dimensional) quotient module by ${\cal F}_0(n|m)=
{\cal F}(n|m)/{\cal F}_1(n|m)$, and (by abuse of notation)
its vectors are denoted by
\[
|l_1,\ldots,l_{n+m}\ra \hbox{ with } l_1+\cdots +l_{n+m}\leq p.
\]
For $h$ an indeterminate ($q=e^h$), the representation of $U_q[sl(n+1|m)]$
into ${\cal F}_0(n|m)$ is irreducible. It is obvious how to 
identify ${\cal F}_0(n|m)$ with $W_p$.

In order to turn ${\cal F}_0(n|m)$ into a unitary $U_q[sl(n+1|m)]$
module, we introduce the Holstein-Primakoff realization.

\begin{theo}[Holstein-Primakoff realization]
Let $p\in\C$. 
The linear map $\varrho : U_q[sl(n+1|m)] \rightarrow W(n|m)$,
defined on the Jacobson generators by
\begin{eqnarray}
\varrho (H_i)&=&p-(-1)^{\t_i}c_i^+c_i^- -\sum_{j=1}^{n+m}c_j^+c_j^-
=p-(-1)^{\t_i}N_i-N, \nn\\
\varrho (a_i^-)&=&q^{N_1+\l +N_{i-1}}\sqrt{{[N_i+1]
\over{N_i+1}}[p-N]}c_i^-, \nn\\
\varrho (a_i^+)&=&\bar{q}^{N_1+\l +N_{i-1}}
\sqrt{{[N_i]\over{N_i}}[p-N+1]} c_i^+,  \label{hp}
\end{eqnarray}
is a homomorphism of $U_q[sl(n+1|m)]$ into $W(n|m)$.
\label{prophp}
\end{theo}

Let us now also consider the special case that $p$ is a positive
integer. Just as in the previous case, the subspace 
${\cal F}_1(n|m)$ is invariant for the action of 
$U_q[sl(n+1|m)]$ under $\varrho$ when $p\in\Z_+$.
It is clearly invariant under the action of $U_q[sl(n+1|m)]$.
Let us consider the following basis of the (finite
dimensional) quotient module ${\cal F}_0(n|m)$~:
\begin{equation}
{(c_1^+)^{l_1}(c_2^+)^{l_2}\l (c_{n+m}^+)^{l_{n+m}}
\over{\sqrt {l_1!l_2!\l l_{n+m}!}}}|0\ra\equiv
|l_1,l_2,\l ,l_{n+m}),\quad l_1+\l +l_{n+m}\leq p.
\label{76}
\end{equation}
It is easy to verify that the action of every $\varrho(H_i)$ and
$\varrho(a_i^\pm)$ upon $|l_1,\ldots,l_{n+m} )$ is the same
as the action of $H_i$ and $a_i^\pm$ in the representation
on $W_p$ given by (\ref{65})-(\ref{66}), under
the identification 
\[
|l_1,\ldots,l_{n+m} ) \equiv |p; l_1,\ldots,l_{n+m} ).
\]
Therefore, it follows that the representation $\varrho$ 
of $U_q[sl(n+1|m)]$
into ${\cal F}_0(n|m)$ (under the Holstein-Primakoff 
realization with $p\in\Z_+$) is an irreducible unitary module when
\[
q=e^{i\phi}\hbox{ with } -{\pi\over p} <\phi < {\pi\over p}.
\]
%In the rest of this section, we shall assume that this is the case.

{}From~(\ref{23}) and~(\ref{hp}), one obtains the Holstein-Primakoff
realization of the remaining Cartan-Weyl elements of 
$U_q[sl(n+1|m)]$~:
\begin{eqnarray}
\varrho(e_{ji})&=&q^{2\t_j(N_j-1)+N_{j+1}+N_{j+2}+\l +N_{i-1}}
\sqrt{{[N_j][N_i+1]\over{N_j(N_i+1)}}}c_j^+c_i^-, \quad (j<i),\nn\\
\varrho(e_{ji})&=&{\bar {q}}^{2\t_iN_i+N_{i+1}+N_{i+2}+\l +N_{j-1}}
\sqrt{{[N_j][N_i+1]\over{N_j(N_i+1)}}}c_j^+c_i^-, \quad (j>i).
\label{77}
\end{eqnarray}

Observe that there is an alternative description of the 
Holstein-Primakoff realization, in terms of deformed 
Bose~\cite{Macfarlane,Biedenharn,Sun,Polychronakos} and
Fermi~\cite{Floreanini2} CAO's $\tilde c_i^\pm$ defined as
\begin{equation}
\tilde{c}_i^-={\sqrt{[N_i+1]\over{N_i+1}}}c_i^-, \quad
\tilde{c}_i^+={\sqrt{[N_i]\over{N_i}}}c_i^+, \quad
\tilde{N}_i=N_i, \quad i\in [1;n+m].
\label{78}
\end{equation}
These elements of $W(n|m)$ satisfy the relations
\begin{equation}
\lb  \tilde{c}_i^-, \tilde{c}_j^+\rb _{q^{\delta_{ij}}}=
\delta_{ij}{\bar q}^{(-1)^{\t_i}\tilde{N}_i}, \quad
[ \tilde{N}_i , \tilde{c}_j^\pm ]=\pm \delta_{ij}
\tilde{c}_j^{\pm}, \quad
\lb  \tilde{c}_i^{\pm}, \tilde{c}_j^{\pm}\rb =
[\tilde{N}_i, \tilde{N}_j ]=0 .
\label{79}
\end{equation}
The Holstein-Primakoff realization can be rewritten 
in terms of these deformed Bose and Fermi operators 
$\tilde{c}_i^{\pm}$. We give it here
for all Cartan-Weyl elements~:
\begin{eqnarray}
\varrho (H_i)&=&p-(-1)^{\t_i}\tilde N_i-\tilde N, \nn\\
\varrho (a_i^-)&=&q^{\tilde{N}_1+\l +\tilde{N}_{i-1}}
\sqrt{[p-\tilde{N}]}\tilde{c}_i^-, \nn\\
\varrho (a_i^+)&=&\bar{q}^{\tilde{N}_1+\l +\tilde{N}_{i-1}}
\sqrt{[p-\tilde{N}+1]} \tilde{c}_i^+, \nn\\
\varrho(e_{ji})&=&q^{2\t_j(\tilde{N}_j-1)+\tilde{N}_{j+1}+
\tilde{N}_{j+2}+\l +\tilde{N}_{i-1}}
\tilde{c}_j^+\tilde{c}_i^-, \quad (j<i), \nn\\
\varrho(e_{ji})&=&{\bar {q}}^{2\t_i\tilde{N}_i+\tilde{N}_{i+1}+
\tilde{N}_{i+2}+\l +\tilde{N}_{j-1}}
\tilde{c}_j^+\tilde{c}_i^-, \quad (j>i). \label{varrho}
\end{eqnarray}
Furthermore, this is easy to extend to a Holstein-Primakoff
realization of $U_q[gl(n+1|m)]$ by
\begin{equation}
\varrho (e_{00})=p-\tilde{N}, \qquad 
\varrho (e_{ii})=\tilde{N}_i.
\end{equation}
The Holstein-Primakoff realization has given us a realization
in terms of oscillators (in~(\ref{hp}) and~(\ref{77}))
or $q$-oscillators (in~(\ref{varrho})).
Observe that this oscillator realization is different from the
one given by Floreanini {\em et al}~\cite{Floreanini2}~: 
in~\cite{Floreanini2} only the Chevalley generators are realized in terms
of oscillators or $q$-oscillators. Furthermore all generators
are bilinear expressions in the oscillators, whereas here the
JGs are linear expressions in the oscillators.

\section{Conclusions}

We have constructed a class of representations of the quantum
superalgebra $U_q[sl(n+1|m)]$, which was also extended to 
$U_q[gl(n+1|m)]$.
Our approach is entirely along the lines of Fock representations
of parastatistics of order $p$, for which the defining relations
are given by~(\ref{3-1})-(\ref{3-3}).
The analogy with parastatistics goes further~: within the Fock
representations, the JGs $a_i^\pm$ can be interpreted as operators
creating or annihilating (quasi)particles, or excitations of
a new kind of quantum statistics.

In order to be more concrete consider a Hamiltonian 
$H=\sum_{i=1}^{n+m}\ \varepsilon_i e_{ii} $. Then, see~(\ref{82b}),
\[
H|p;r_1,\ldots,r_{n+m})
=\sum_{i=1}^{n+m}\ \varepsilon_i r_i
|p;r_1,\ldots,r_{n+m}). 
\]
Therefore the vector $|p;r_1,\ldots,r_{n+m})$ can be interpreted as a
state consisting of $r_1$ particles
with energy $\varepsilon_1$, $r_2$ particles
with energy $\varepsilon_2$, and so on, $r_{n+m}$ particles
with energy $\varepsilon_{n+m}$. Moreover, according to
(\ref{65})-(\ref{66}) any operator $a_i^+$ (resp. $a_i^-$) creates
(resp. annihilates) a particle on the orbital~$i$. Since
$r_i\in \Z_+ $ for $i\in[1;n]$ and
$r_{i}\in \{0,1\}$ for $i \in[n+1;m]$
the particles on the first $n$ orbitals behave like bosons,
and the particles on the next orbitals like fermions.
This is however not quite the case if $p<n+m$, since 
$\sum_{i=1}^{n+m}r_i\le p.$
In other words the system cannot accommodate more than $p$
particles. Therefore the statistics falls in the group
of exclusion statistics in the broad sense~\cite{Haldane}~: the number
of particles to be accommodated on a certain orbital 
depends on the number of particles already accommodated in the system.
What are the properties of the underlying statistics is a question still
to be answered.

Another property worth to be studied is to analyze
what happens if $p \rightarrow \infty$ and $q \rightarrow 1$. 
Having in mind the results from~\cite{PV00} we expect that in this limit
the operators $A(p,q)_i^\pm=a_i^\pm/{\sqrt p}$ become genuine
Bose CAO's for $i\in [1;n]$ and genuine Fermi CAO's for 
$i\in [n+1;n+m]$. If so, then for large $p$ and 
values of $q$ close to 1 the operators $A(p,q)_i^\pm$ 
describe small deviations from the canonical quantum statistics. Moreover 
these CAO's are defined in a state space with positive definite scalar 
product. Among the various noncanonical statistics (Gentile intermediate 
statistics~\cite{Gentile}, parastatistics~\cite{Green}, 
infinite statistics~\cite{Haag}, parons~\cite{Gr},
quons~\cite{ArikGr}) only the quons have the same property.
Therefore parallel to quons the $A(p,q)_i^\pm$ operators
may appear as another  candidate to describe eventual
small violations of canonical quantum statistics 
(see~\cite{Gr} where also experiments for detecting small violations of
statistics are discussed). 

We believe also (having in mind again the results in~\cite{PV00})
that the CAO's of $U_q[sl(n+1|m)]$ and their Fock
representations will be natural ``building blocks"
for any multicomponent $t-J$ supersymmetric lattice model.
To this end we note that at each site $i$
the Hubbard operators $X^{0k}$ and $X^{k0}$~\cite{Fo89} 
(we suppress the site index) are nothing but nondeformed Jacobson 
generators $a_k^-$ and $a_k^+$, respectively.
Then the representations with $p=1$  satisfy the hard-core 
condition forbidding configurations with two or more particles to be
accommodated simultaneously on each lattice site. 

Some of the results related to this quantum statistics were
already published in an earlier paper~\cite{PSV00}.
Let us underline the new contributions in the present paper.
Theorem~1 (section~II) was already stated without proof in~\cite{PSV00}, 
since it is the main ingredient to describe the quantum statistics;
here we have given its relevant background and a complete proof.
Sections~III and~IV contain our key results; all of them are original.
We have constructed a class of representations of $U_q[sl(n+1|m)]$
labelled by an arbitrary number $p$.
When $p$ is a positive integer, this representation
is indecomposable and the corresponding quotient module is finite
dimensional.
The derivation of the action of the JGs on basis elements of
these representations is highly nontrivial.
In section~IV we have selected the unitary representations,
with respect to the (in physics) natural Hermiticity 
condition~(\ref{dagger}) considered as an anti-involution,
and the requirement that the usual Fock space metric should
be positive definite.
It is interesting to note that the selected representations
remain irreducible when $q$ is a root of unity.

The Dyson and Holstein-Primakoff realizations of $U_q[sl(n+1|m)]$
were given in an earlier paper by one of us~\cite{Palev99}, but only for
the Chevalley generators.
Here, in section~V, we give the realization for all Cartan-Weyl
elements of $U_q[sl(n+1|m)]$. Such realizations are relevant since
also in the classical case ($q=1$) the realization of all Cartan-Weyl
generators (Bargmann-Schwinger realizations, ladder representations)
are of physical importance.
Observe that it is far from trivial to deduce the realization
for all Cartan-Weyl elements from those of the realization for the
Chevalley generators. This would be rather hard because the
expressions of all Cartan-Weyl elements in terms of the Chevalley
generators are very involved and difficult to manage, see
e.g.~\cite{PalevTolstoy}.
In the present case, the problem was overcome because we were able
to give the Dyson and Holstein-Primakoff realizations of the Jacobson
generators of $U_q[sl(n+1|m)]$. 
Since the expressions of the remaining Cartan-Weyl elements
in terms of the Jacobson generators is simple,
the Dyson and Holstein-Primakoff realizations of all Cartan-Weyl
elements followed without too much trouble.

\section*{Acknowledgements}
T.D.\ Palev and N.I.\ Stoilova were supported by NATO (Collaborative
Linkage Grant) during their visit to Ghent.
T.D.\ Palev also wishes to acknowledge Ghent University for a
visitors grant.
N.I.\ Stoilova is thankful to Prof.\ H.D.\ Doebner for constructive
discussions and to the Humboldt Foundation for its support.

\end{document}